\documentclass[letterpaper,10pt,conference]{ieeeconf}
\IEEEoverridecommandlockouts
\usepackage{float}
\usepackage{graphicx}
\usepackage{amssymb}

\usepackage{amsmath}
\usepackage{mathtools}  
\mathtoolsset{showonlyrefs}
\usepackage{algorithmic}
\usepackage{algorithm}
\usepackage{nccmath}

\usepackage{amsthm}
\theoremstyle{plain}

\begin{document}

\title{\LARGE \bf Self-Organization in Networks: A Data-Driven Koopman  Approach}

\author{Claudia Caro-Ruiz,
		Duvan Tellez-Castro,
        Andr\'es Pavas,
        and Eduardo Mojica-Nava
\thanks{Claudia Caro-Ruiz, Duvan Tellez-Castro, Andr\'es Pavas, and Eduardo Mojica-Nava are with Departamento de Ingenier\'ia El\'ectrica, Universidad Nacional de Colombia, Bogot\'a, Colombia 111321, $\{$\texttt{clccaroru, datellezc, fapavasm, eamojican$\}$@unal.edu.co}}        
\thanks{Duvan Tellez-Castro is supported in part by Colciencias 727-2016. Claudia Caro-Ruiz is supported in part by Colciencias 647-2015. }
}        
 
\maketitle

%%%%%%%%%%%%%%%%%%%%%%%%%%%%%%%%%%%%%%%
%%%%%%%%%%%%   Abstract   %%%%%%%%%%%%%
%%%%%%%%%%%%%%%%%%%%%%%%%%%%%%%%%%%%%%%
\begin{abstract}

 %%%%%%%%%%%%%%%%%%%%%%%%%%%%%%%%%%%%%%%%%%
%%%%%%%%%%%%   Abstract %%%%%%%%%%%%%%%%%%
%%%%%%%%%%%%%%%%%%%%%%%%%%%%%%%%%%%%%%%%%%
Networks out of equilibrium present embedded dynamics that are characterized by multiple equilibria and sudden transitions.  Global changes in dynamics behavior arise by strong coupling between node neighbors producing that each node leaves its natural stable state and joins to an organized global activation. Significative changes in dynamics occur when network approach a regime transition. 
In this paper, we present a data-driven approach based on Koopman analysis spectrum for the identification of local patterns and dynamical properties near regime shifts. To illustrate these ideas we propose two applications. First, an Integral-and-Fire (IFO) coupled oscillator network showing transitions due to synchronization. Second, a Bak-Sneppen model of evolution that depicts the network patterns tending to self-organized criticality (SOC). We obtain Koopman modes and spectra by a dynamic mode decomposition techniques applied to data measurements. By this, we identify local self-organizing properties and transitions indicators from the dynamical properties of the proposed applications.
%%%%%%%%%%%%%%%%%%%%%%%%%%%%%%%%%%%%%%%%%%
 
\end{abstract}
%%%%%%%%%%%%%%%%%%%%%%%%%%%%%%%%%%%%%%%
%%%%%%%%%%%%   KeyWords %%%%%%%%%%%%%%%
%%%%%%%%%%%%%%%%%%%%%%%%%%%%%%%%%%%%%%%
%\IEEEpeerreviewmaketitle
%%%%%%%%%%%%%%%%%%%%%%%%%%%%%%%%%%%%%%%
%%%%%%%%%%%   INTRODUCTION %%%%%%%%%%%%
%%%%%%%%%%%%%%%%%%%%%%%%%%%%%%%%%%%%%%%
%%%%%%%%%%%%%%%%%%%%%%%%%%%%%%%%%%%%%%%%%%
%%%%%%%%%%%   INTRODUCTION %%%%%%%%%%%%%%%
%%%%%%%%%%%%%%%%%%%%%%%%%%%%%%%%%%%%%%%%%%
\section{Introduction}
\label{sec:intro}
%%%%%%%%%%%%%%%%%%%%%%%%%%%%%%%%%%%%%%%%%%
Critical transitions in networks are a major indicator of complexity in system dynamics as is described by  \cite{dorogovtsev2008critical}. Any system dynamics approaching to its critical region shows a strong sensitivity to external perturbations and parameter variations of its micro-scale dynamics. By this, stability, and predictability of system dynamics can change significantly.  Control and operation actions designed for systems out of critical regime cannot manage these changes correctly. 

Self-organization is an implicit mechanism governing the occurrence of transitions in complex networks.
Even if the dynamics are tending to a stable and cooperative regime, or if it is tending to criticality (see \cite{noel2014bottom},  \cite{wang2016growth}), the dynamics in networks are subjected to the emergence of local patterns of behavior that changes in a discontinuous way.  Results from \cite{asllani2014theory}, present an extension of the theory of pattern formations for a reaction-diffusion system over directed networks. Also, in  \cite{nicolaides2016self},  a study of quantized self-organized network patterns over a Swift-Hohenberg continuum model is used to depict complex suite localized patterns. 
Also, there are some other spatial patterns that can arise before a critical transition. For instance, scale-invariance distributions of avalanche clusters occurrence,  the general trend towards increasing spatial coherence, and  the increasing of nodes cross-correlation as is presented in \cite{scheffer2009early}. 
Also,  \cite{moon2015network} presents a way to rebuild network properties by using spatial patterns and coherency of network near a critical transition. Moreover, the identification of changes in local patterns can be used to predict the occurrence of transitions  as is described by \cite{zhang2015predictability}. 
Even the results described before  are developed on very specific models, the identification of these patterns and shifts in real large scale complex systems cannot be made analytically. This is because dimension and complexity evade the modeling and identification process. By approach this problem the use of data-driven techniques would be an alternative for the analysis and identification of these patterns.  
The Koopman operator is a linear and infinite-dimensional and acting over observables of dynamic systems. Applications and theoretical contributions from Koopman analysis can be observed for example in stability analysis, power system analysis,   and control design \cite{Mauroy2014}, \cite{Susuki2016},  \cite{Brunton2016}.  %\cite{Susuki2011a} describe stability assessment of %power systems without models, a monitoring  power system dynamics using dynamic mode decomposition is presented by \cite{mohapatra2016fast}, and a power grid partitioning is presented by \cite{raak2014partitioning}.\\

In this paper, we present a data-driven approach based on Koopman analysis spectrum for the identification of local patterns and dynamical properties near regime shifts. We applied the technique to networks tending to synchronization dynamics and to SOC. By this, we identify local self-organizing properties and transitions indicators  of networks that can be obtained from data measures.  
The rest of the paper is organized as follows. Section \ref{sec:Net_Phenomena} presents the characteristics of the  studied phenomena. Section \ref{sec:koopman} presents the Koopman approach. Section \ref{sec:results} presents the numerical results of the proposed method. Finally, Section \ref{sec:conclusions} presents some conclusions and future work.

%%%%%%%%%%%%%%%%%%%%%%%%%%%%%%%%%%%%%%%
%%%%%% Review of the literature %%%%%%%
%%%%%%%%%%%%%%%%%%%%%%%%%%%%%%%%%%%%%%%
%%%%%%%%%%%%%%%%%%%%%%%%%%%%%%%%%%%%%%%%%%
%%%%%%%%%%%   Review of the literature %%%%%%%%%%%%%%%
%%%%%%%%%%%%%%%%%%%%%%%%%%%%%%%%%%%%%%%%%%
\section{Network Phenomena}
\label{sec:Net_Phenomena}
%%%%%%%%%%%%%%%%%%%%%%%%%%%%%%%%%%%%%%%%%%
In this section, we describe two models that depict self-organization in their dynamics. Both models are used to obtain data-measurements useful for the analysis of self-organization patterns and regime shifts in networks. 

\subsection{Synchronization Patterns: Integral-and-Fire Coupled Oscillators}
Considers a laticce network of size $N\times N$ where each node $i$  is characterized by two states variables: a state variable $\theta_i$ that increases linearly with time, and an energy variable $E_i$,   depending on $\theta$ and the node network interactions. The model is described by \cite{corral1995self}. The variables dynamics are described by,
\begin{equation}\label{}
\frac{d\theta_i}{dt}=1,
\end{equation}
\begin{equation}
\frac{dE_i}{dt}= \gamma\left(K-E_i\right),
\end{equation}
where $K=1/\left(1-e^{\left(-\gamma\right)}\right)$ and $E_i\left(\theta_i\right)=K\left[1-e^{\left(\gamma\theta_i\right)}\right]$
to achieve $E_i\left(\theta_i=0\right)=0$ and $E_i\left(\theta_i =1 \right)= 1$.

Once a node becomes critical achieving ($E_{i}\geq E_c = 1$) it fires and transfer energy to its neighbors according to the following rule
\begin{align*}
E_{j}&\rightarrow E_{j}+\epsilon,\\
E_{i}&\rightarrow 0,
\end{align*}
where $\epsilon$ is the strength of the coupling and $j \in \mathcal{N}_i$ and $\mathcal{N}_i$ is the neighborhood of node $i$. By this, some of the node neighbors may become critical and generate an avalanche that propagates through the lattice. The local dynamics are interrupted once an avalanche starts; also, when the avalanche finishes, the  local driving actions works again. This behavior presents a time scale separation.  One scale for the slow dynamics of the internal phase state variables and the other for the instantaneous interaction between the units (pulse coupling). 
In the model, there exist an intrinsic dynamics leading the network to the threshold. When a node reaches the threshold, it fires producing a change in the phase of its neighbors. It could produce new firings generating an avalanche.  Here, we study the state of the system after the avalanches, not the avalanche size distribution as is usual in Self-organized criticality nevertheless the ingredients are the same. Under this conditions, SOC and synchronization are very similar.

These IFO oscillators had been used to study the behavior of cells, neurons, and other biological systems. Conditions as nonlinear convex driving described by $\gamma$ for the own dynamics and long-range interactions between nodes are conditions that warranty a stationary state that presents complete synchronization. Moreover under particular driven and interaction state it can also give a self-organized criticality behavior.

\subsection{Self-Organized Criticality: Bak-Sneppen Model }
The term Self-Organized Criticality depicts the dynamics of a system tending to a critical point. Thus, the macroscale behavior of the network displays an spatial and temporal scale-invariance characteristic. It means that a power law relation between avalanche size and frequency indicates that the stationary state of the attractor is critical.  SOC behavior is present in many examples in nature as landscapes, forest fires, geodesic formations, and earthquakes.  Here we use a model presented in literature called the Bak-Sneppen Model \cite{bak1993punctuated}.

 The Bak-Sneppen model is a simple mathematical model of biological macroevolution. It describes the adaptation process of interacting species. The entire model evolves to a self-organized criticality state where periods of non-evolution alternate with avalanches of extinctions producing evolutionary changes.  Extinctions of all sizes, including mass extinctions, may be a simple consequence of ecosystem dynamics.  The models work based on two simple rules. First, by finding the species with the lowest fitness and randomly changing its fitness. Second, at the same time changing the fitness of the species on the lower species neighborhood. By the time where a threshold minimum fitness is achieved, it shows the highest critical steady state where the macro dynamics of the system has evolved.  Also, avalanches of fitness changes occur below the threshold.
 The setup for the Bak-Sneppen model is a ring of size $N$, that represent ecological niches, each occupied by a species population. Associated with each site is a fitness $x_i\in U[0,1]$. The network fitness evolves according to the following rules:
 \begin{enumerate}
 \item Find node index $i_{min}$ associate to
 \begin{equation}
x_{i_{min}}=\min\left\{x_i|i=1,...,N\right\}.
 \end{equation}
 \item Replace $x_{i_{min}}$ and its nearest neighbors by random numbers taken from $U[0,1]$
 \item{Go back to step 1}
 \end{enumerate}
 Thus, the system starts with a uniform fitness distribution $\rho(0)=U[0,1]$, then applying the algorithm generates consecutive distributions $\rho(1),\rho(2)$ until achieving a critical distribution $\rho(\infty)=U[x_{crit},1]$ where $\rho(\infty)<1$.

%%%%%%%%%%%%%%%%%%%%%%%%%%%%%%%%%%%%%%%
%%%%%%%%%%%%%%% Methods %%%%%%%%%%%%%%%
%%%%%%%%%%%%%%%%%%%%%%%%%%%%%%%%%%%%%%%
%%%%%%%%%%%%%%%%%%%%%%%%%%%%%%%%%%%%%%%%%%
%%%%%%%%%%%%%%%   Methods %%%%%%%%%%%%%%%%
%%%%%%%%%%%%%%%%%%%%%%%%%%%%%%%%%%%%%%%%%%
\section{Koopman Analysis}
\label{sec:koopman}
%%%%%%%%%%%%%%%%%%%%%%%%%%%%%%%%%%%%%%%%%%
%%%%%%%%%%%%%%%%%%%%%%%%%%%%%%%%%%
%%%%%%%%%%%SEC-Koopman%%%%%%%%%%%%
%%%%%%%%%%%%%%%%%%%%%%%%%%%%%%%%%%
In this section, we definite briefly the Koopman operator,  its basic spectral properties, and also we present the Koopman mode decomposition.

\subsection{Koopman Operator}
 Consider a dynamical system  as generating of a set $\mathcal{X}$ and self-map $F:\mathcal{X}\rightarrow \mathcal{X}$. For instance,
\begin{align}
	\dot{\bold{x}}=F(\bold{x}), \hspace*{1 cm} \bold{x} \in \mathcal{X}, \label{eq1} 
\end{align}
  where  $\mathcal{X}$ is a nonempty compact space and the map $F$ is continuous, and let $\varphi(\bold{x}_0,t)$ be the flow  that represents the trajectories of \eqref{eq1} with initial condition $\bold{x}_0$. Now, consider a family of functions called observables $g:\mathcal{X}\rightarrow \mathbb{C}$ where $g \in \mathbb{G}$ is infinite dimensional Hilbert space and $\mathbb{C}$ is the complex set. A particular observable for systems with a stable fixed point is a \textit{ function of  Lyapunov}. Rather than a space state transformation, we consider its Koopman operator $\mathcal{K}:=\mathcal{K}_{F}$ defined by
\begin{align}
	\mathcal{K}_{F}g:=g \circ F,
\end{align}
for some observables on  $\mathcal{X}$.

\begin{align}
	\mathcal{K}g=g(\varphi(\bold{x}_0,t)),\hspace*{1 cm} \mathcal{K}g \in \mathcal{X}.\label{eq2}
\end{align}

A compact form of  Equation \eqref{eq2} can be represented by
\begin{align}
	\dot{g}=\mathcal{K}g, \hspace*{1 cm} \mathcal{K}g \in \mathcal{X}.\label{eq3} 
\end{align}
 
 In other words, the Koopman operator acts in an infinite dimensional linear operator that describes the evolution of observables along of trajectories in the space $\mathcal{X}$, unlike $F(\bold{x})$ acting in the finite dimensional nonlinear space states $\mathbb{R}^n$  \cite{Budisic2012}, \cite{Mauroy2013a}.

This operator captures full information for a large class of nonlinear dynamical systems. A comparison between traditional and Koopman analysis of dynamical systems is shown in Table \ref{tab:comparison_table}.
\begin{table}[!h] 
\caption{Koopman analysis vs Traditional  approachs.}
\label{tab:comparison_table}
\begin{center}
\begin{tabular}{ll}
\hline%----------------------------------------------------------------- 
\textbf{Koopman analysis}        &  \textbf{Traditional methods}      \\
\hline%----------------------------------------------------------------- 
Linear dynamics and infinite     &  Nonlinear dynamics and finite     \\
dimensional eigenfunction space. &  dimensional state space.          \\
\hline%----------------------------------------------------------------- 
Full spectral analysis.          & Linearization techniques.  \\
\hline%-----------------------------------------------------------------  
Free model, operators is         & Identification methods and         \\
approximated or analyzed from    & model techniques.                  \\
 experimental simulation data.    &                                    \\
\hline %----------------------------------------------------------------- 
Hard connection with             & Directly relation with             \\
physical intuition.              & physics system.                    \\
 \hline%-----------------------------------------------------------------                           
Analysis and synthesis is a      & Extensive theory.                  \\
challenging problem.             &                                    \\
\hline%-----------------------------------------------------------------    
 \end{tabular}
\end{center}
\end{table}

\subsection{Spectral Analysis}
We begin by redefining \eqref{eq1} in a system in time discrete  in order to reconstruct the Koopman spectrum from a set of data, through of a map the state $\bold{x_0}$ forward time $t$ into the future to $x(t_0 + t)$ following, 
\begin{align}
	\varphi(\bold{x}_0,t)=\bold{x}(t_0 + t)=\bold{x}_0+\int_{t_0}^{t_0 + t}F(\bold{x}(\tau ))d \tau. \label{eq5}
\end{align}

Then, we obtain a new discrete-time dynamical system given by
 \begin{align}
 \bold{x}_{k+1}=	\varphi (\bold{x}_k). \label{eq6}
 \end{align}

The dynamical system defined in \eqref{eq6} and the one defined by $\mathcal{K}$ in \eqref{eq3}  are analogous, the connection between theses is the identity operator $\bold{x}=\bold{g(x)}$ and the Koopman spectrum given by $(\mu_k,\phi_k,v_k)$, which are Koopman eigenvalues, eigenfunctions, 
and modes respectively.  We begin by defining a basis for space $\mathcal{X}$, and represent our observable functions $\bold{g(x)}$ in this basis as is described by \cite{Tu2013}.
\begin{align}
	\bold{g(x)}=\sum^\infty_{k=1}v_k\phi_k(\bold{x}). 
\end{align} 

In doing this, is possible obtain the evolving of systems \eqref{eq1} as
\begin{align}
	F(\bold{x})=(\mathcal{K}g)(\bold{x})=\sum^\infty_{k=1}v_k(\mathcal{K}\phi_k)(\bold{x})=\sum^\infty_{k=1}e^{\mu_k \Delta_t} v_k \phi_k(\bold{x}), \label{eq4}
\end{align}
where $\Delta_t$ is a fixed sample time \cite{Williams2015}.

The Koopman eigenvalues $\mu_k$ described of temporal attributes of $v_k$, i.e.,  fast or slow dynamics are determined by the real part of $\mu_k$ and the imaginary part represented  its frequency, or  magnitude and phase in continuous case \cite{Rowley2009}. In addition,   each $v_k$ is presented as the projection onto the eigenfunction  space and in many cases is assumed that the system has an attractor, in another way each $v_k$ oscillates with a single frequency as is formulated by \cite{Susuki2016} and \cite{Mezic2012}. 
 
%Our goal is to analyze the dynamics of \eqref{eq1} using the spectrum Koopman for the identification of self-organization threshold in different kind of systems. 

\begin{figure}[!t]
\centering
\includegraphics[scale=0.43]{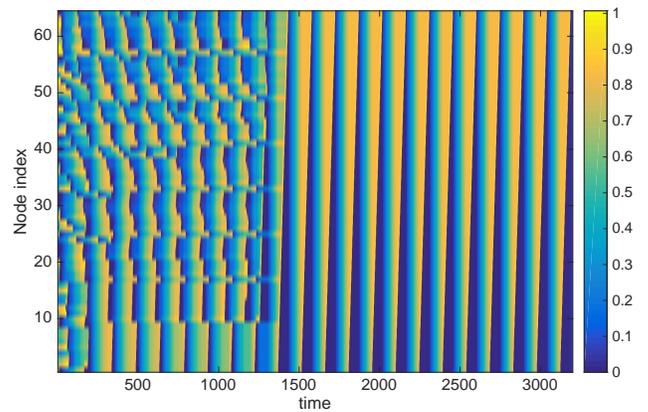}
\caption{Dynamics of  IFO coupled laticce network of size N=64. The node index  in vertical axis indicates  phase state variable  for each node in time.  During the 1000 to 1200 iterations, a regime transition occurs.  IFO parameters used for obtained the data are $\epsilon =0.145$, $\gamma=2$, and iteration step $\Delta t =0.01$s.}\label{Sync}
\end{figure}

 \begin{figure}[!t]
\centering
\includegraphics[scale=0.43]{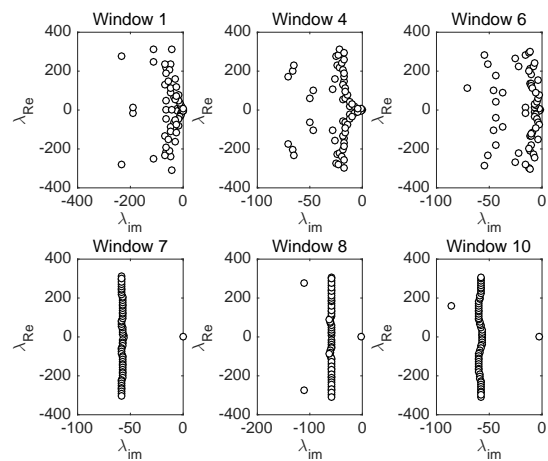}
\caption{Koopman eigenvalues for the IFO network. Time windows from data measures are taken every 200-time iterations. The figure depicts transition on spectral patterns between the eighth and the ninth time window.}\label{KooplambdaSYnc}
\end{figure}

 \begin{figure}[!t]
\centering
\includegraphics[scale=0.43]{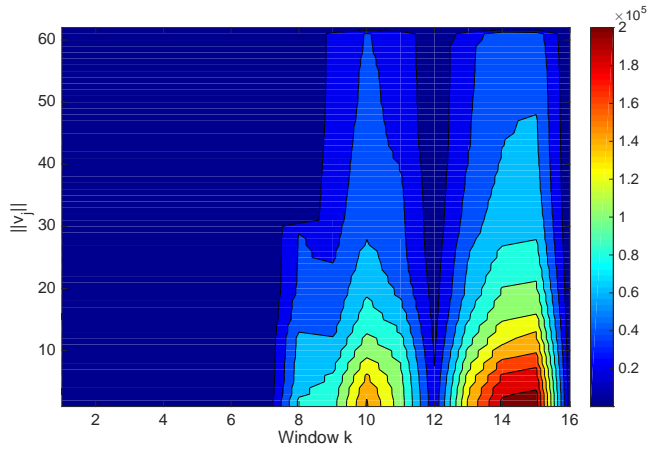}
\caption{Koopman modes normalized amplitude. Transition regime is detectable by the modes magnitudes. An increasing of five orders of magnitude occurs during the synchronization transition. After the transitions, amplitude maintains higher, considering that modes size before the transition is lower than unity.}\label{KoopmodesSync}
\end{figure}

\begin{figure}[!t]
\centering
\includegraphics[scale=0.43]{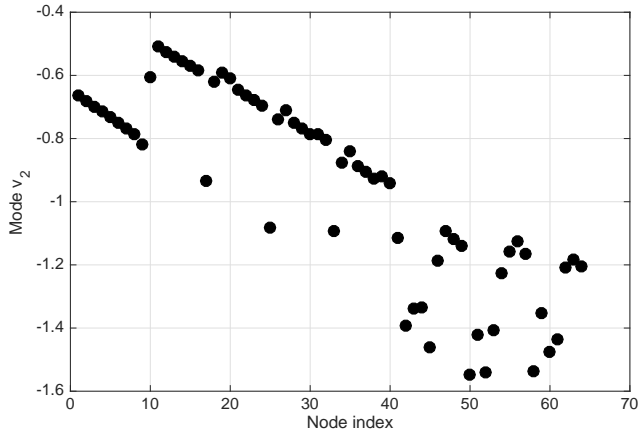}
\caption{Koopman dominant mode $v_2$. Near the threshold, mode two shows patterns of interaction and properties of the network structure.}\label{Koopmanmode2}
\end{figure}

\subsection{Computation of Koopman Modes}
Computation of Koopman  spectrum is a challenging problem. The algorithms more usually are the Dynamic Mode Decomposition (DMD) and Extended Dynamic Mode Decomposition. Both uses a set of data for  obtaining eigenvalues and eigenvectors of linear operators but the second includes a dictionary for finding eigenfunctions allowing us to find  the full spectrum Koopman, for details see \cite{Tu2013} and \cite{Williams2015}.

%%%%%%%%%%%%%%%%%%%%%%%%%%%%%%%%%%%%%%%
%%%%%%%%%%%%%% Results %%%%%%%%%%%%%%%%
%%%%%%%%%%%%%%%%%%%%%%%%%%%%%%%%%%%%%%%
%%%%%%%%%%%%%%%%%%%%%%%%%%%%%%%%%%%%%%%%%%
%%%%%%%%%%%%%   Results %%%%%%%%%%%%%%%%%%
%%%%%%%%%%%%%%%%%%%%%%%%%%%%%%%%%%%%%%%%%%
\section{Simulation Results}
\label{sec:results}
%%%%%%%%%%%%%%%%%%%%%%%%%%%%%%%%%%%%%%%%%%
\subsection{Numerical Results for IFO Coupled Network}

Fig. \ref{Sync} presents the dynamic behavior of a set of coupled IFO oscillators sharing energy over a lattice network of size 64.  Without any external driving, the system displays in the stable state relaxation oscillations formed by large assemblies of units all with the same phase. In this way, the system achieves complete synchronization.  Avalanche patterns occurring in different places from the network follows self-organizing patterns until a point of transition where coherence emerges. In the first stage of the graphic, it shows various groups of elements subject to avalanches of different sizes. The shift occurs after 1200 time iterations. Fig. \ref{KooplambdaSYnc} presents results from applying Koopman spectrum for the analysis of self-organization in this example of synchronization. The figure presents nine-time windows where each time window $[t_0,t_k]$ presents a stage of the IFO coupled dynamics evolving in time.   By using DMD,  we obtain Koopman spectrum with 62 Koopman modes that decompose the entire dynamics of the system.  Fig. \ref{KooplambdaSYnc}   shows the obtained Koopman eigenvalues. Time windows from 1 to 6 present dynamics before the relaxing oscillations transition. 
It is possible to observe for each stage that there exist different pattern dynamics depending on the nodes activated by the avalanches. The eigenvalues grouped near the imaginary axes  are continuously changing its configuration.  During the eight-time window, a sudden change in the network behavior occurs and the oscillators achieve the synchrony. In this way, an avalanche occupying the entire systems occurs periodically in the system.  Fig. \ref{KooplambdaSYnc}  presents the change in Koopman eigenvalues behavior during the transition and after that. In the synchronizing state, the eigenvalues organize with the same decaying rate. Many eigenvalues are far away from the imaginary axes. The system achieves a stable fixed point defined by the Koopman spectrum depicted in the figure. The governing dynamics is split  into two kinds described by eigenvalues.  A group with fast time scale dynamic characterized by the highest decaying rate and a  group with slow time scale dynamic from the period of the synchronized oscillations represented by the eigenvalue near the imaginary axis.  The fast time scale represents the constant size avalanches presented in each network oscillation. 
 \begin{figure}[!t]
\centering
\includegraphics[scale=0.43]{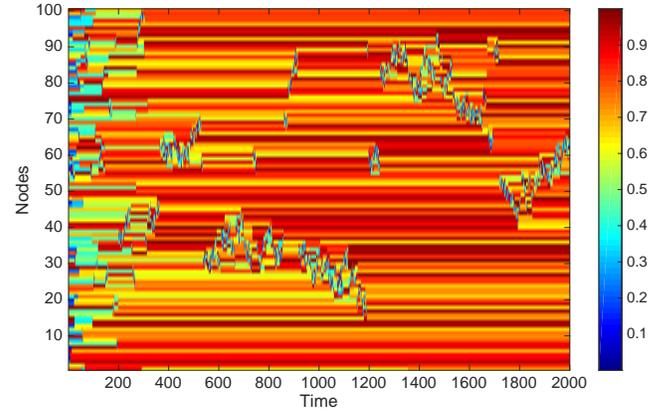}
\caption{Population dynamics for a Bak-Sneppen self-organized criticality model over a ring network of size $N=100$. Avalanches of different sizes occur in time. During time iterations 1000 and 1200, the system achieves the  self-organization threshold.}\label{BKlatice}
\end{figure}

 Using the Koopman spectrum we can obtain an analysis related to Koopman dominant modes. For each Koopman eigenvalue, there exist and associate Koopman mode.  Modes magnitude can be organized and normalized to identify participation of dominant modes in the network dynamics. In Fig.  \ref{KoopmodesSync}   we depict the behavior of Koopman modes for data measures obtained from IFO coupled oscillators during twenty-five-time windows.  For data taken after the critical transition, Koopman modes participation is very slow. Its magnitude values are smaller than unity in many cases.  During the time window before the synchronization transition occurs, a sudden abrupt increase of four orders of magnitudes for modes magnitude occurs. This sudden change in half of the modes is a signal of arriving the synchronization threshold in the network. After that, the Koopman modes amplitude continues being of this order of magnitude. Also, the participation of every mode is similar and significative, presenting a strong correlation between dynamics of every network node.  Also, the dominant second mode associated to real eigenvalue displays information about structure and patterns of network dynamics.  Fig.  \ref{Koopmanmode2} presents the Koopman dominant mode two $\|v_2\|$ obtained from the dynamic mode decomposition applied to data in the time frame before the transition. Its components are related to any node of the network.  Patterns of interaction between nodes are implicit in this mode. Before the synchronization shift, relaxing oscillations of small size are occurring over the network.  Nodes index $i \in [1,9]$ have mode components grouped sharing an avalanche pattern. In a similar way, nodes with index $i \in [11,40]$ present a repetitive local pattern of interaction. Also nodes from $i\in [4,64]$, they share similar regional patterns inside the former groups. 
\begin{figure}[!t]
\centering
\includegraphics[scale=0.4]{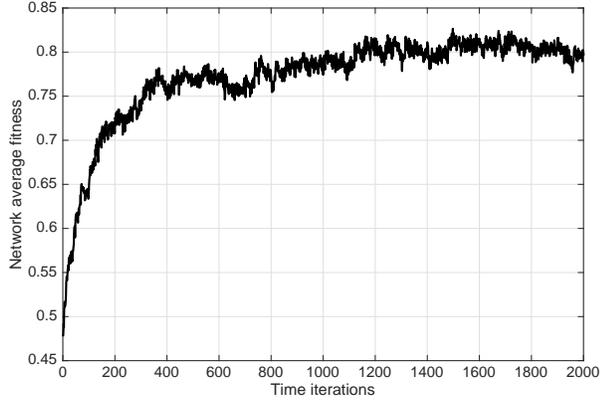}
\caption{Average fitness for the Bak-Sneppen self-organized criticality model over a ring network of size $N=100$. Average fitness achieves its maximum during the transition regime and maintains the same value along the time. }\label{fitness}
\end{figure}
\begin{figure}[!t]
\centering
\includegraphics[scale=0.43]{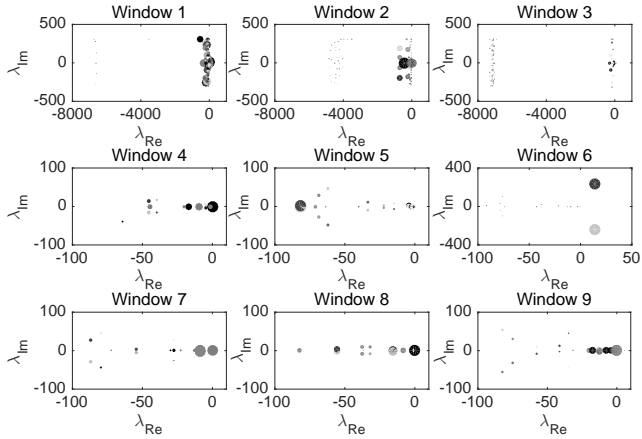}
\caption{Koopman spectrum for data snapshots of the Bak-Sneppen model. Unstable regimes are intermittent in the spectrum. It depicts a sudden transition for the sixth time window. }\label{BKeigenvalues}
\end{figure}
In this section, we present the results obtained from the application of Koopman analysis approach to data measurements obtained from networks dynamics showing synchronization and criticality thresholds. 
 
\subsection{Numerical Results for Bak-Sneppen Self-Organized Criticality Model}
Fig. \ref{BKlatice} presents the dynamics of the Bak-Sneppen model over a network ring with a population of 100 species. The figure shows the changes in  nodes fitness due to random local changes. The network fitness achieves an average threshold of 0.8 approximately. In Fig. \ref{fitness}, we observe the network average fitness in time. It has a smooth behavior in comparison with the random sudden changes that occurs at the local level.  In this system, iteration between species emerges through a kind of red queen effect that optimizes the entire system behavior. We analyze the data measurements obtained for the Bak-Sneppen model by the use of the Koopman spectrum. 

\begin{figure}[!t]
\centering
\includegraphics[scale=0.43]{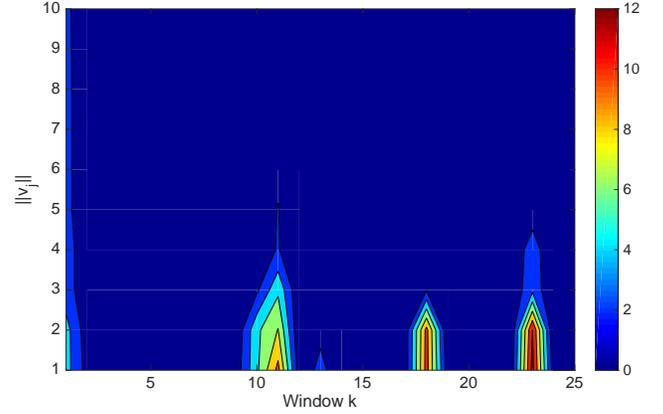}
\caption{Dominant mode amplitudes for twenty-five data snapshots. Higher values in dominant modes magnitudes occur after the transition. During the shift, dominant modes take values near to zero.}\label{BKmodescountour}
\end{figure}

We use the DMD numerical method to get spectrum and modes of Koopman operator. Fig. \ref{BKeigenvalues} presents the Koopman spectrum taken for time windows of two hundred iterations. Each window displays a different behavior in the spectrum governed by two groups of eigenvalues; large decaying rates eigenvalues that represent fast dynamics in avalanches, and small increasing rate eigenvalues related to the growing behavior in fitness dynamics. The critical point of self-organization can be related to spectrum behavior during the sixth time window. The size of markers indicators for each eigenvalue represents the magnitude of the corresponding modes. Growing rates govern the system behavior during the transition generating instability and bigger avalanches. It is the same time window where the system achieves the maximum average fitness for the network.  In the same way, Fig. \ref{BKmodescountour} presents an oscillating behavior in modes amplitude.  Koopman dominant modes are few for each iteration. Also, its magnitudes change significantly depending on the kind of dynamics that are governing the behavior.  If growing rates govern the dynamics, modes dominant amplitudes tend to zero. Stable modes driven by high decaying rates magnitudes represents a dynamic of large size avalanches, and modes dominated by small decaying rates accounts for a network governed by constant average fitness and lower size avalanches. 

\begin{figure}[!t]
\centering
\includegraphics[scale=0.43]{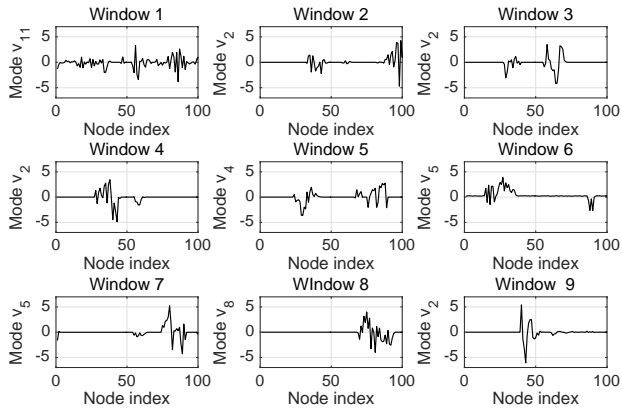}
\caption{Modes components shows the zones where avalanches occur during a time lapse.}\label{BKmodesavalanches}
\end{figure}

Fig. \ref{BKmodesavalanches} shows an useful result obtained from the observations of Koopman modes. The figure displays the governing zero-frequency mode component for each time window. These modes reveal the local self-organized patterns occurring during the avalanches. The mode component value for each node quantifies its participation on local interacting patterns during the avalanches in a time frame. This information could be used to rebuild structural properties of the network and to identify nodes centrality regarding stable and unstable modes.  
%%%%%%%%%%%%%%%%%%%%%%%%%%%%%%%%%%%%%%%
%%%%%%%%%%% Discussion %%%%%%%%%%%%%%%%
%%%%%%%%%%%%%%%%%%%%%%%%%%%%%%%%%%%%%%%
%\input{6_Discussion} 
%%%%%%%%%%%%%%%%%%%%%%%%%%%%%%%%%%%%%%%
%%%%%%%%%%% CONCLUSIONS %%%%%%%%%%%%%%%
%%%%%%%%%%%%%%%%%%%%%%%%%%%%%%%%%%%%%%%
%%%%%%%%%%%%%%%%%%%%%%%%%%%%%%%%%%%%%%%%%%
%%%%%%%%%%%%%   Conclusions %%%%%%%%%%%%%%%
%%%%%%%%%%%%%%%%%%%%%%%%%%%%%%%%%%%%%%%%%%
\section{Conclusions}
\label{sec:conclusions}
%%%%%%%%%%%%%%%%%%%%%%%%%%%%%%%%%%%%%%%%%%
In summary, we describe IFO coupled oscillators and Bak-Sneppen model as useful applications for the study of self-organization patterns and dynamical regime shifts in the network. We present a data-driven analysis of synchronization and SOC thresholds for the described applications. For this analysis, we use a Koopman Operator approach. From a DMD method, we obtain the Koopman spectrum and modes, and we use them to analyze dynamical properties concerning to self-organization. 
The main results obtained in this work are the identification of a stable fixed point for synchronization dynamics described by Koopman eigenvalues, and the use of Koopman dominant modes amplitude as a signal for prediction of transitions. Also, we identify regime changes for SOC models by the utilization of the Koopman operator and the usefulness of principal mode components for the identification of local patterns generated from avalanches and node participation to local self-organized events. 
As future works, we will describe by the use of Koopman operator theory the attractors related to regime transitions and its relation with network properties. Also,  we will use Koopman operator approach for the design of control systems to reduce size and frequency of avalanches in SOC phenomena.

%%%%%%%%%%%%%%%%%%%%%%%%%%%%%%%%%%%%%%%
%%%%%%%%%%% Acknowledgment %%%%%%%%%%%%
%%%%%%%%%%%%%%%%%%%%%%%%%%%%%%%%%%%%%%%
%\section*{Acknowledgment}
%The authors would like to thank Prof ... for  letting us ... and his help that greatly improved the manuscript.
%%%%%%%%%%%%%%%%%%%%%%%%%%%%%%%%%%%%%%%
%%%%%%%%%%%%% References %%%%%%%%%%%%%%
%%%%%%%%%%%%%%%%%%%%%%%%%%%%%%%%%%%%%%%
\bibliography{KOOPSOC, koopmanbiblio}

\begin{thebibliography}{10}
\providecommand{\url}[1]{#1}
\csname url@samestyle\endcsname
\providecommand{\newblock}{\relax}
\providecommand{\bibinfo}[2]{#2}
\providecommand{\BIBentrySTDinterwordspacing}{\spaceskip=0pt\relax}
\providecommand{\BIBentryALTinterwordstretchfactor}{4}
\providecommand{\BIBentryALTinterwordspacing}{\spaceskip=\fontdimen2\font plus
\BIBentryALTinterwordstretchfactor\fontdimen3\font minus
  \fontdimen4\font\relax}
\providecommand{\BIBforeignlanguage}[2]{{%
\expandafter\ifx\csname l@#1\endcsname\relax
\typeout{** WARNING: IEEEtran.bst: No hyphenation pattern has been}%
\typeout{** loaded for the language `#1'. Using the pattern for}%
\typeout{** the default language instead.}%
\else
\language=\csname l@#1\endcsname
\fi
#2}}
\providecommand{\BIBdecl}{\relax}
\BIBdecl

\bibitem{dorogovtsev2008critical}
S.~N. Dorogovtsev, A.~V. Goltsev, and J.~F. Mendes, ``Critical phenomena in
  complex networks,'' \emph{Reviews of Modern Physics}, vol.~80, no.~4, p.
  1275, 2008.

\bibitem{noel2014bottom}
P.-A. No{\"e}l, C.~D. Brummitt, and R.~M. D'Souza, ``Bottom-up model of
  self-organized criticality on networks,'' \emph{Physical Review E}, vol.~89,
  no.~1, p. 012807, 2014.

\bibitem{wang2016growth}
Y.~Wang, H.~Fan, W.~Lin, Y.-C. Lai, and X.~Wang, ``Growth, collapse, and
  self-organized criticality in complex networks,'' \emph{Scientific reports},
  vol.~6, 2016.

\bibitem{asllani2014theory}
M.~Asllani, J.~D. Challenger, F.~S. Pavone, L.~Sacconi, and D.~Fanelli, ``The
  theory of pattern formation on directed networks,'' \emph{Nature
  communications}, vol.~5, 2014.

\bibitem{nicolaides2016self}
C.~Nicolaides, R.~Juanes, and L.~Cueto-Felgueroso, ``Self-organization of
  network dynamics into local quantized states,'' \emph{Scientific reports},
  vol.~6, 2016.

\bibitem{scheffer2009early}
M.~Scheffer, J.~Bascompte, W.~A. Brock, V.~Brovkin, S.~R. Carpenter, V.~Dakos,
  H.~Held, E.~H. Van~Nes, M.~Rietkerk, and G.~Sugihara, ``Early-warning signals
  for critical transitions,'' \emph{Nature}, vol. 461, no. 7260, pp. 53--59,
  2009.

\bibitem{moon2015network}
H.~Moon and T.-C. Lu, ``Network catastrophe: Self-organized patterns reveal
  both the instability and the structure of complex networks,''
  \emph{Scientific reports}, vol.~5, 2015.

\bibitem{zhang2015predictability}
X.~Zhang, C.~Kuehn, and S.~Hallerberg, ``Predictability of critical
  transitions,'' \emph{Physical Review E}, vol.~92, no.~5, p. 052905, 2015.

\bibitem{Mauroy2014}
A.~Mauroy and I.~Mezi{\'{c}}, ``{Global stability analysis using the
  eigenfunctions of the Koopman operator},'' \emph{arXiv.org}, vol. 9286,
  no.~c, pp. 1--34, 2014.

\bibitem{Susuki2016}
Y.~Susuki, I.~Mezic, F.~Raak, and T.~Hikihara, ``{Applied Koopman operator
  theory for power systems technology},'' \emph{Nonlinear Theory and Its
  Applications, IEICE}, vol.~7, no.~4, pp. 430--459, 2016.

\bibitem{Brunton2016}
S.~L. Brunton, B.~W. Brunton, J.~L. Proctor, and J.~N. Kutz, ``{Koopman
  invariant subspaces and finite linear representations of nonlinear dynamical
  systems for control},'' \emph{PLoS ONE}, vol.~11, no.~2, 2016.

\bibitem{corral1995self}
{\'A}.~Corral, C.~J. P{\'e}rez, A.~D{\'\i}az-Guilera, and A.~Arenas,
  ``Self-organized criticality and synchronization in a lattice model of
  integrate-and-fire oscillators,'' \emph{Physical review letters}, vol.~74,
  no.~1, p. 118, 1995.

\bibitem{bak1993punctuated}
P.~Bak and K.~Sneppen, ``Punctuated equilibrium and criticality in a simple
  model of evolution,'' \emph{Physical review letters}, vol.~71, no.~24, p.
  4083, 1993.

\bibitem{Budisic2012}
M.~Budi{\v{s}}i{\'{c}}, R.~Mohr, and I.~Mezi{\'{c}}, ``{Applied Koopmanism},''
  \emph{Chaos}, vol.~22, no.~4, 2012.

\bibitem{Mauroy2013a}
A.~Mauroy and I.~Mezi{\'{c}}, ``{A spectral operator-theoretic framework for
  global stability},'' \emph{Proceedings of the IEEE Conference on Decision and
  Control}, no.~3, pp. 5234--5239, 2013.

\bibitem{Tu2013}
J.~H. Tu, C.~W. Rowley, D.~M. Luchtenburg, S.~L. Brunton, and J.~N. Kutz, ``{On
  dynamic mode decomposition - theory and applications},'' \emph{Journal of
  computational dynamics}, no. September, pp. 1--30, 2013.

\bibitem{Williams2015}
M.~O. Williams, I.~G. Kevrekidis, and C.~W. Rowley, ``{A Data???Driven
  Approximation of the Koopman Operator: Extending Dynamic Mode
  Decomposition},'' \emph{Journal of Nonlinear Science}, vol.~25, no.~6, pp.
  1307--1346, 2015.

\bibitem{Rowley2009}
C.~W. Rowley, I.~Mezi{\'{c}}, S.~Bagheri, P.~Schlatter, and D.~S. Henningson,
  ``{Spectral analysis of nonlinear flows},'' \emph{Journal of Fluid
  Mechanics}, vol. 641, pp. 115--127, 2009.

\bibitem{Mezic2012}
I.~Mezi{\'{c}}, ``{Analysis of Fluid Flows via Spectral Properties of the
  Koopman Operator},'' \emph{Annual Review of Fluid Mechanics}, vol.~45, no.~1,
  p. 121005161233001, 2012.

\end{thebibliography}
%%%%%%%%%%%%%%%%%%%%%%%%%%%%%%%%%%%%%%%
\end{document}